# Increasing the Peneteration of Renewables by Releasing Merchant Energy Storage Flexibility

*Farhad Samadi Gazijahani*
Department of Electrical Engineering
Azarbaijan Shahid Madani University
Tabriz, Iran
f.samadi@azaruniv.ac.ir

*Rasoul Esmaeilzadeh*
Department of Distributed Generation
Azarbaijan Regional Electric Company
Tabriz, Iran
rasoul_zadeh@yahoo.com

*Abstract*—It is anticipated that large-scale energy storage facility (ESF) will become as an essential part of the future energy markets to increase the penetration of renewables. In this paper, a new optimization algorithm is developed to participate large-scale *merchant* ESF to deal with renewables fluctuations. The ESF is an *investor-owned* and *independently-operated*, seeking to maximize its total profit, while hedging the system against net-load variation. To do so, an exact computationally efficient bi-level mixed integer linear programming (MILP) is envisaged. The ESF operator tries to maximize its own profit at upper level while ISO seeks to clear the markets at lower level. In principle, the ESF will be able to construct its bidding/offering curves to buy/sell from/to different markets. To maintain the tractability of the problem, each lower level is replaced by its KKT optimality conditions and nonlinearities are converted into linear equivalents using strong duality theory rendering the single-level MILP and consequently recast as an MPEC. Additionally, the proposed framework is constructed according to the information gap decision theory (IGDT) instrument to capture the adverse impact of uncertainty on the profit of ESF.

*Keywords—Flexibility, Energy storage; Renewable; Transmission network, Optimization.*

## I. INTRODUCTION

It is anticipated that a large portion of future electricity markets will be supplied by renewable resources as *priority areas* all around the world [1]. In this regard, energy storage facilities (ESFs) can be regarded as one of the best solutions to store the excess renewable energy when is not needed, and give back to the system in emergency conditions ensuring more reliable power supply and facilitating the penetration of renewable resources [2]-[3]. As a result, significant benefits will directly be achieved by employing ESFs, for example, reducing operation costs, flatting electricity consumption, decreasing renewable curtailments, and providing ancillary services for systems [4].

Renewable energy, typically photovoltaic power and wind power, will have higher percentages in total power generation. With the increasing percentages of wind power and photovoltaic power, the impacts of extreme weather events on them would be critical to improve the entire system resilience. For example, photovoltaic power could be a means to supply local loads when the main networks are in outages due to hurricanes. However, hurricanes indicate severe cloud cover, which directly reduces the photovoltaic power outputs.

On the other hand, the evolution of energy markets has become more dependent on renewables due to climate change and environmental issues caused by utilizing the fossil fuels, as well as the depletion and cost of conventional power plants [5]. Toward this direction, ESFs are being deployed not only to arbitrage in the electricity markets but also to curb the uncertainties of renewables [6]. Locations and sizes of distributed generators could have great impacts on system resilience improvement. Two-stage models can be used to deal with this issue.

Technically speaking, the ESF can be utilized in two different ways so called commercial and technical ESF [7]. Commercial ESF has generally been used to earn economic interests through participation in the different power markets while technical ESF demodulates the technical limitations of the existing network to achieve these goals. In fact, in the economic analysis the impact of the network topology does not take into account regardless of technical ESF that the real time influence of the local network considers to investigate the effect of the ESF on the technical specifications of the system [8].

In some works, ESFs have been examined to enable the provision of system services in the power markets such as reliability increasing, providing reserve, and frequency control [9]-[11]. Moreover, in some others, the ESF has been used for economic targets through participation in energy and reserve markets [12]-[15]. Besides, in existing literature, some methods have been developed to incorporate the uncertainty induced by renewables into the problem by scenario-based programming [16], point estimate method [17], robust optimization [18], and chance constrained approach [19].

However, multiple works have been conducted on the ESF participation in the power markets, but there are so many shortcomings that need to be addressed properly. The present work aims to address the remaining problems of prior literature by providing the following contributions:

1) A bi-level MILP is propounded for ESF to participate in ramp-product market aimed at minimizing renewable swings.
2) For constituting a computationally tractable model, KKT optimality conditions and strong duality theory are used to recast the problem as an MILP by applying MPEC.
3) An information gap decision theory (IGDT) has been extended to curb the uncertainties induced by prices.

The remainder of this paper is categorized as follows. Section 2 demonstrates the modelling of ESF; the problem formulation is presented in Section 3; the concept of IGDT approach is explained in Section 4 and the results are shown and discussed in Section 5. Finally, the relevant conclusions are provided in Section 6.

## II. ESF Modelling

Energy storage is also an important means to increase system resilience against different natural disasters, and determining locations and sizes for fixed storages is one of critical point. The price-maker ESF operator purchases energy at low-price hours and then sells it to the market at high-price periods. Under these situations, ESF operator strives to maximize its own profit by arbitrage in both markets. In this regard, this paper examines a novel bi-level MILP to make the proposed bidding problem computationally tractable along with implementation of IGDT to mitigate the risk of participation in the power markets. The main feature of IGDT compared with other uncertainty modeling approaches is that IGDT is composed of two immunity functions that considers both positive and negative aspects of uncertainty.

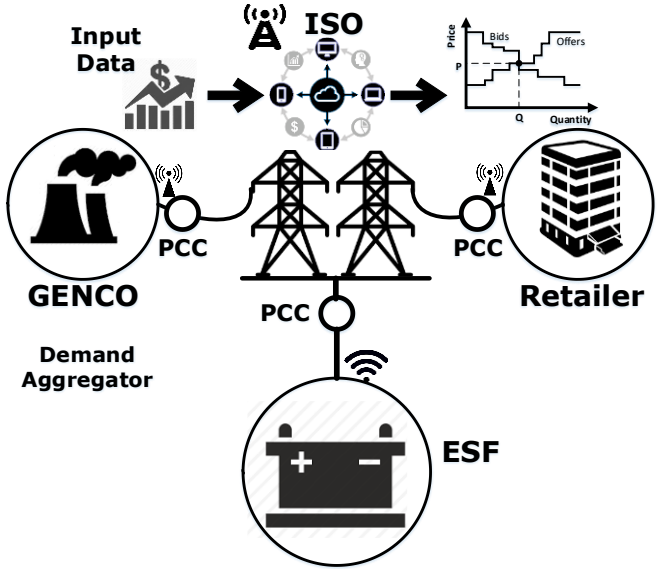

Fig. 1. The correlated ESF aggregated for participating into the market.

The state of charge (SOC) of ESF is shown by equation (1). In addition, the allowable bound of SOC and charging/discharging constraints of ESF can be restricted by (2)-(4). The initial SOC of ESF is given in (5) and the power of BES determines through (6). Moreover, constraint (7) is employed to eschew charging and discharging of EES at the same time. Finally, equation (8) calculates the operation cost of EES over scheduling horizon.

$$SOC_t^{ES} = SOC_{t-1}^{ES} + \Delta t \left( P_{t,ch}^{ES} \cdot \eta_{ch}^{ES} - P_{t,dis}^{ES} / \eta_{dis}^{ES} \right) \quad (1)$$

$$SOC_{\min}^{ES} \leq SOC_t^{ES} \leq SOC_{\max}^{ES}, \quad \forall t \quad (2)$$

$$P_{ch,\min}^{ES} \leq P_{t,ch}^{ES} \leq P_{ch,\max}^{ES}, \quad \forall t \quad (3)$$

$$P_{dis,\min}^{ES} \leq P_{t,dis}^{ES} \leq P_{dis,\max}^{ES}, \quad \forall t \quad (4)$$

$$SOC^0 = SOC^T \quad (5)$$

$$P_t^{ES} = P_{t,ch}^{ES} - P_{t,dis}^{ES}, \quad \forall t \quad (6)$$

$$x_{t,ch}^{ES} + y_{t,dis}^{ES} \leq 1, \quad \forall t \quad (7)$$

$$K^{ES} = C^{ES} \times \left( x_{t,ch}^{ES} \cdot P_{t,ch}^{ES} + y_{t,dis}^{ES} \cdot P_{t,dis}^{ES} \right), \quad \forall t \quad (8)$$

Note that the merchant ESF is a price maker facility that is capable of altering price by exercising its market power, i.e., physical and economic withholding, to profitably raise market prices for sustained periods of time. For this testimony, the price quota curves have been used [20] to link the market price proportion to the generation and demand quotas of price maker ESF which has intuitively been illustrated in Fig. 2.

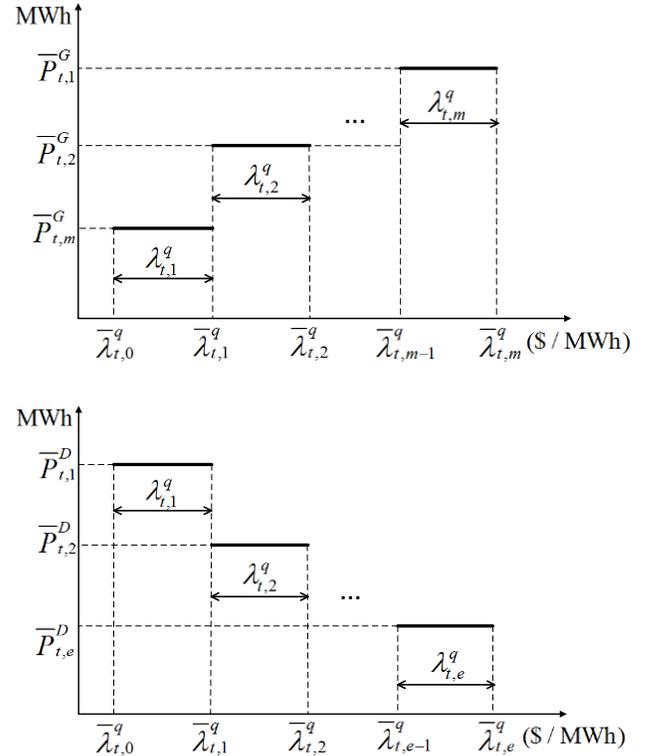

Fig. 2. Generation and demand price-quota curves of price-maker ESF.

It is assumed that the pay-as-bid pricing scheme is employed to settle both oligopoly markets. Unlike the uniform pricing, in a pay-as-bid auction instead of paying market price to winning players, they will receive the price which have offered. Therefore, in pay-as-bid auctions, the players instead of offering their marginal costs, they offer their own bidding curves based on their individual forecasts of the price in order to increase their own profit as much as possible [21]. That is why, the pay-as-bid auction are referred to as "discriminatory auctions" since they pay winners different prices in comparison to market price. It should be notified that in several European markets like in Germany and Italy, the balancing markets are settled under a pay-as-bid pricing scheme [22]. Note that a change in market clearing mechanism will have small effect on other market concerns, consisting ensuring adequate transmission networks and increasing the role of ESF.

## III. PROBLEM FORMULATION

This paper deals with optimal offering strategy of a merchant ESF, who aims to participate in both DA and RT oligopoly markets as a price maker agent. The main objective of the ESF is to maximize its own profit. For this purpose, a new bi-level optimization model is formulated, where the uncertainty lies in the market prices. The proposed bi-level framework consists of an upper level that maximizes ESF profit and two lower levels that represent the DA and RT market-clearing process.

### 1.1. Formulation of Upper Level

The main objective of the upper level is to maximize the profit of price maker ESF by arbitrage in both DA and RT oligopoly markets as below:

$$\underset{profit}{Max}\ F_{Up}^{ES} = \sum_t \left\{ \left( O_t^{DA,ES}.P_t^{DA,ES} + O_t^{RT,ES}.P_t^{RT,ES} \right) - \left( B_t^{DA,ES}.D_t^{DA,ES} + B_t^{RT,ES}.D_t^{RT,ES} \right) \right\} \quad (9)$$

Equation (9) consists of four major parts, where two first terms denote revenues obtained from participation in DA and RT oligopoly markets, i.e., offering strategy, and the next two terms also show the purchasing power required by ESF from DA and RT markets, i.e., bidding strategy. It should be emphasized that the proposed model is linear and convex because DA and RT market prices determine in the lower levels. At the end, constraints (10)-(13) restrict the offering/bidding power of ESF in the DA and RT markets, respectively.

$$P_t^{DA,ES} \geq 0,\ \forall t:\ O_t^{DA,ES} \quad (10)$$

$$P_t^{RT,ES} \geq 0,\ \forall t: O_t^{RT,ES} \quad (11)$$

$$D_t^{DA,ES} \geq 0,\ \forall t: B_t^{DA,ES} \quad (12)$$

$$D_t^{RT,ES} \geq 0,\ \forall t: B_t^{RT,ES} \quad (13)$$

### 1.2. Formulation of Lower Levels

The market clearing process is run at the lower-level for both DA and RT markets on the basis of social welfare maximization. Axiomatically, all players who intend to participate in the market submit a list of bidding power blocks and their corresponding prices to the ISO for every hour of the planning horizon. Then, the ISO clears the markets be fulfilling a single-round or multi-round unit commitment problem considering security constraints of the system. In this paper, the DA market clearing (14) consists of four major parts: two first terms show the incomes achieved by players participated in the DA market, and also the next two parts express the bidding costs of players to purchase the power from DA market.

$$Min\ F_{L1}^{DA} = \sum_t \begin{Bmatrix} O_t^{DA,ES}.P_t^{DA,ES} + \sum_{j \in GENCO} O_{t,j}^{DA}.P_{t,j}^{DA} - \\ B_t^{DA,ES}.D_t^{DA,ES} - \sum_{i \in Retailer} B_{t,i}^{DA}.D_{t,i}^{DA} \end{Bmatrix} \quad (14)$$

$$D_t^{DA,ES} + \sum_i D_{t,i}^{DA} - P_t^{DA,ES} - \sum_j P_{t,j}^{DA} = 0,\ \forall t:\lambda^{DA} \quad (15)$$

$$0 \leq P_t^{DA,ES} \leq P_t^{\max},\ \forall t:\ \delta_t^{\min}, \delta_t^{\max} \quad (16)$$

$$0 \leq D_t^{DA,ES} \leq D_t^{\max},\ \forall t: \kappa_t^{\min}, \kappa_t^{\max} \quad (17)$$

$$0 \leq P_{t,j}^{DA} \leq P_j^{\max},\ \forall j:\ \mu_j^{\min}, \mu_j^{\max} \quad (18)$$

$$0 \leq D_{t,i}^{DA} \leq D_i^{\max},\ \forall i:\ \xi_i^{\min}, \xi_i^{\max} \quad (19)$$

$$-f_l^{\max} \leq f_l \leq f_l^{\max},\ \forall l, \forall t:\ \omega_i^{\min}, \omega_i^{\max} \quad (20)$$

Equation (15) guarantees the power mismatch in the DA market. The constraints (16)-(17) show the limitations on the selling and purchasing powers by ESF and also, constraints (18)-(19) restrict the permissible selling/purchasing power for GENCO and retailers, respectively. Finally, constraint (20) enforces allowable bound of power transactions in the lines of transmission network. The dual variables are declared at the corresponding constraints following a colon. The primal variables of DA lower level problem are $\Xi_{DA}^{Primal} = \{D_t^{DA,ES}, P_t^{DA,ES}, D_{t,i}^{DA}, P_{t,j}^{DA}\}$ and its dual variables are $\Xi_{DA}^{Dual} = \{\lambda^{DA}, \delta_t^{\min}, \delta_t^{\max}, \kappa_t^{\min}, \kappa_t^{\max}, \mu_j^{\min}, \mu_j^{\max}, \xi_i^{\min}, \xi_i^{\max}, \omega_i^{\min}, \omega_i^{\max}\}$ Similar to the DA market clearing formulation, the objective function of RT market clearing can be applied.

Regarding to bi-level structure, the first-order KKT optimality conditions is replaced instead of lower-level problems resulting in single level MPEC problem. To this end, primal-dual transformation technique has been executed by inserting a collection of primal and dual constraints in exchange of initial equations. Equations (21)-(31) denote the dual constraints of DA market clearing problem while (32) is its strong duality equality that specifies the equality of primal and dual variables at the optimal point.

$$O_t^{DA,ES} - \lambda^{DA} + \delta_t^{\max} - \delta_t^{\min} = 0,\ : P_t^{DA,ES},\ \forall t \quad (21)$$

$$O_{t,j}^{DA} - \lambda^{DA} + \mu_j^{\max} - \mu_j^{\min} = 0,\ : P_{t,j}^{DA},\ \forall t \quad (22)$$

$$-B_t^{DA,ES} + \lambda^{DA} + \kappa_t^{\max} - \kappa_t^{\min} = 0,\ \forall t: D_t^{DA,ES} \quad (23)$$

$$-B_{t,i}^{DA} + \lambda^{DA} + \xi_i^{\max} - \xi_i^{\min} = 0: D_{t,i}^{DA},\ \forall i \quad (24)$$

$$\lambda_l^s - \lambda_l^r - \omega_l + \omega_l^{\max} - \omega_l^{\min} = 0: f_l,\ \forall l \quad (25)$$

$$\delta_t^{\max} \geq 0;\ \delta_t^{\min} \geq 0\ \ : \overset{\max}{\delta_t}, \overset{\min}{\delta_t}\ \forall t \quad (27)$$

$$\mu_j^{\max} \geq 0;\ \mu_j^{\min} \geq 0\ \ : \overset{\max}{\mu_j}, \overset{\min}{\mu_j}\ \forall j \quad (28)$$

$$\kappa_t^{\max} \geq 0;\ \kappa_t^{\min} \geq 0\ \ : \overset{\max}{\kappa_t}, \overset{\min}{\kappa_t}\ \forall t \quad (29)$$

$$\xi_i^{\max} \geq 0;\ \xi_i^{\min} \geq 0\ \ : \overset{\max}{\xi_i}, \overset{\min}{\xi_i}\ \forall i \quad (30)$$

$$\omega_l^{\max} \geq 0;\ \omega_l^{\min} \geq 0\ \ : \overset{\max}{\omega_l}, \overset{\min}{\omega_l}\ \forall l \quad (31)$$

$$\begin{aligned} & O_t^{DA,ES} P_t^{DA,ES} + \sum_j O_{t,j}^{DA}.P_{t,j}^{DA} - B_t^{DA,ES} D_t^{DA,ES} \\ & -\sum_i B_{t,i}^{DA}.D_{t,i}^{DA} + P_t^{\max}.\delta_t^{\max} + P_j^{\max}.\mu_j^{\max} + \\ & D_t^{\max}.\kappa_t^{\max} + D_i^{\max}.\xi_i^{\max} + f_l^{\max}.\omega_l^{\max} = 0: \Delta_t^{DA} \end{aligned} \quad (32)$$

Constraints (33)-(42) illustrate the complementarity of KKT conditions corresponding to the inequalities of MPEC of the DA market clearing mechanism.

$$0 \leq P_t^{DA,ES} \perp \delta_t^{\min} \geq 0, \quad \forall t \tag{33}$$

$$0 \leq \left(P_t^{\max} - P_t^{DA,ES}\right) \perp \delta_t^{\max} \geq 0, \quad \forall t \tag{34}$$

$$0 \leq P_{t,j}^{DA} \perp \mu_t^{\min} \geq 0, \quad \forall t, \forall j \tag{35}$$

$$0 \leq \left(P_j^{\max} - P_{t,j}^{DA}\right) \perp \mu_t^{\max} \geq 0, \quad \forall t, \forall j \tag{36}$$

$$0 \leq D_t^{DA,ES} \perp \kappa_t^{\min} \geq 0, \quad \forall t \tag{37}$$

$$0 \leq \left(D_t^{\max} - D_t^{DA,ES}\right) \perp \kappa_t^{\max} \geq 0, \quad \forall t \tag{38}$$

$$0 \leq D_{t,i}^{DA} \perp \xi_t^{\min} \geq 0, \quad \forall t, \forall i \tag{39}$$

$$0 \leq \left(D_i^{\max} - D_{t,i}^{DA}\right) \perp \xi_t^{\max} \geq 0, \quad \forall t, \forall i \tag{40}$$

$$0 \leq f_l \perp \omega_l^{\min} \geq 0, \quad \forall t, \forall l \tag{41}$$

$$0 \leq \left(f_l^{\max} - f_l\right) \perp \omega_l^{\max} \geq 0, \quad \forall t, \forall l \tag{42}$$

Note that strong duality theory and big-M method have been deployed to linearize the non-linear terms of the model [23]. Finally, a single level MILP model is reached by replacing the KKT optimality conditions of lower levels into the upper level. Fig. 3 outlines an overview of the proposed scheme.

It is worth noting that the proactive actions include pre-positioning of energy storage units and DG resources, as well as implementing demand response to deal with upcoming events and net-load variations. Multiple studies have shown that the deployment of mobile energy storage units and DGs at network stations in the pre-disturbance phase reduces load shedding and speeds up network recovery in the post-disturbance phase. DG is also one of the most effective approaches to enhance the resilience of distribution systems. It should be noted that DG reroutes power flow by changing the status of manual and automatic switches, thus leading to enhanced system resilience. Therefore, in the model proposed in this paper, proactive actions including pre-positioning of mobile energy storage units and DG are performed to enhance the system flexibility during emergency conditions.

## IV. IGDT TECHNIQUE

Due to the unpredictability of solar and wind power and voltage-current conditions in the power network, there may be a reduction in the network reliability and a direct impact on the electricity market. To minimize the negative consequences of new participants, greater control and automation of the distribution network are needed, which enable faster decision-making. Ultimately, with greater DER integration, there is a greater requirement for adequate planning and operation. Therefore, it will be necessary to take advantage of large data flows to ensure a reliable and efficient network, considering its unpredictability. To plan the operation of modern distribution networks, TSO will have to use not only consumption forecasts but also production forecasts of new participants, and then forecasts of available flexibility services.

In the pool markets, because of several uncertain resources such as unpredictable offering strategies of rivals and sharp fluctuation of renewables, the decision-making will be treated as a complex process affecting by multiple factors [24]. Choosing an appropriate method among existing methods for uncertainty modeling should be done consciously so that being compatible with existing information about the nature of uncertainty and also must justify the required implementing time. For instance, in the stochastic approaches such as scenario-based programming the probability density function of uncertain parameters should be known which their executing is very time consuming or in the fuzzy modelling methods like α-cut technique the membership function of uncertain parameters should be given which requires running multiple simulations for different levels of membership degrees [25]. To overcome the shortcomings of conventional uncertainty modelling approaches, a novel IGDT based approach has been introduced to tackle the severe uncertainty of the problem [26].

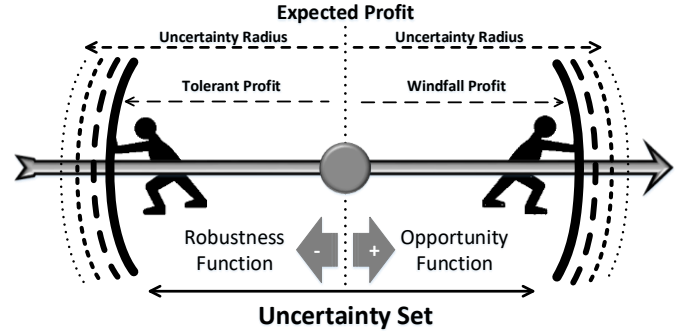

Fig. 4. Concepts of IGDT method based on robustness/opportunity functions.

The IGDT intends to discover the decision variables so that maximizes the chance of success as much as possible [27]. It should be emphasized that the IGDT can be regarded in the form of risk-averse strategy, i.e., robustness function, and opportunity seeker strategy, i.e., opportunity function. In the robustness function the main objective is to set lower profit target that can be accepted by maximizing the tolerance against the uncertainty, while the main objective of opportunity function is to reach higher profit target to pursue through augmenting the chance of success. In other words, the *robustness function* represents the maximum level of uncertainty in which the minimum acceptable profit can be obtained (i.e., tolerant profit). On the contrary, the *opportunity function* refers to the minimum level of uncertainty that allows for inclusive success (i.e., maximum windfall profit) [28]. This concept for modelling of severe uncertainties has graphically been illustrated in Fig. 4.

$$Z_{BC} = \max_{\Phi} \left\{ f(\Phi, X) \right\} \tag{43}$$

$$M(\Phi, X) \geq 0 \tag{44}$$

$$H(\Phi, X) = 0 \tag{45}$$

$$X \in \Gamma\left(\overline{X}, \alpha\right) \tag{46}$$

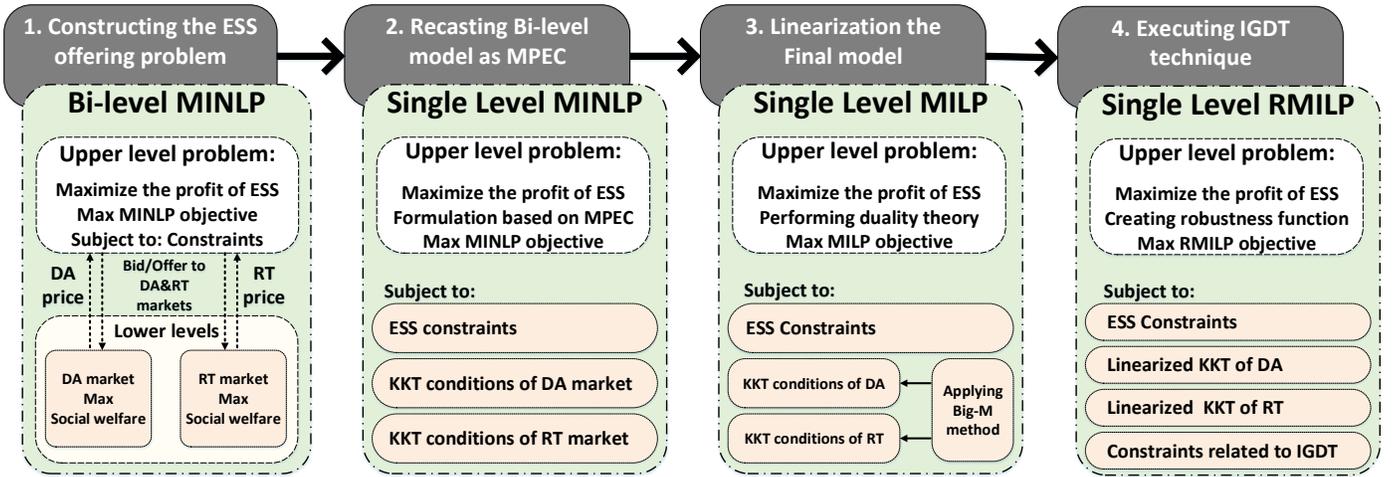

Fig. 3. Comprehensive outline of the model for ESF participation problem.

The initial objective function, i.e., basic case, determines in (43) that should be maximized. Further, (44)-(45), respectively, correspond to the inequality and equality constraints of the problem and (46) mentions decision variables. The robustness and opportunity functions can be formulated as (47)-(48), based on critical profit (49) or permissible profit (50).

$$\alpha_t(Z, f_R) = \max_F \left\{ \alpha : \min_{\mathbb{Z}} F(\Xi, \lambda_t) \geq f_R \right\} \quad (47)$$

$$\beta_t(Z, f_O) = \min_F \left\{ \alpha : \max_{\mathbb{Z}} F(\Xi, \lambda_t) \geq f_O \right\} \quad (48)$$

$$f_R = (1 - \beta_R) Z_{BC} \quad (49)$$

$$f_O = (1 + \beta_O) Z_{BC} \quad (50)$$

$$\lambda_t \in \left\{ \lambda_t^{DA}, \lambda_t^{Bal} \right\} \quad (51)$$

$$\mathbb{Z}(\alpha, \overline{\lambda}_t) = \left\{ \lambda_t : \left| \frac{\lambda_t - \overline{\lambda}_t}{\overline{\lambda}_t} \right| \leq \alpha \cdot \Gamma(t) \right\}, \alpha \geq 0 \quad (52)$$

The vector of uncertain parameters has been given in (51), where their envelope-bound info-gap model has been enforced by (52).

## V. SIMULATION RESULTS

In this section, the proposed model for arbitrage strategy of a price-maker ESF is implemented on the 30-bus transmission network [29] and numerical results are presented here to validate the effectiveness of the model. In particular, different aspects of problem such as ESF power rating and uncertainty modelling by means of IGDT approach are analyzed. Planning and forecasting have always been needed by TSO, and they find their application, among other things, in creating hourly production curves and in planning the purchase and sale of electricity.

The price-maker ESF aims to maximize its own profit by arbitrage in both oligopoly DA and RT markets. The maximum purchasing/selling power for a given hour which ESF operator can purchase/sold from/to both markets is 200/100 MW. At each time, the ESF operator by accomplishing an OPF in ESF, determines that amount of selling/purchasing power to/from markets. In the case of ESF faces with power shortages, the ESF operator participates in the markets as a retailer and submits its bid to the market in order to supply its demand. On the contrary, when there is power surplus in the ESF, the operator of ESF decides to sell it in the market, therefore, it participates as a GENCO and submits its offer to the market. Under this situation, existing an appropriate tool to reduce consumption of ESF at high-price periods seems to be most beneficial.

In addition, the optimal charging/discharging scheduling of ESF facility has been depicted in Fig. 5. On the other side, Table I is shown to clarify the impact of market power on the obtained revenues, with taken into consideration of both price-maker and price-taker strategies. As deduced from this Table, by executing market power via ESF, its revenue will be considerably increased because of altering market price in direction of maximizing its own profit as much as possible. Actually, market power increases the flexibility of ESF in conjunction with its rivals who have participated in the market.

Fig. 4 illustrates the deviations of LMP compared to initial values with (a) and without (b) consideration of ESF. It can be seen that congestion has caused locational price diversity from the initial levels. Specifically, the LMP at bus 2 is relatively lower than the initial price, which is located on one side of the congested line. Besides, the LMP at bus 4 is considerably higher than the initial price, which is located on the other side of the congested line. It should be stated that under congestion, the ESF conducts peak shaving not only *across time* but also to some extent *across location*.

TABLE I
Cost-benefit analysis of different participation model of ESF facility

|  | Modelling | |
|---|---|---|
|  | Price-maker | Price-taker |
| Revenue in DA ($) | 9782.66 | 8165.02 |
| Revenue in RT ($) | 3477.21 | 2508.47 |
| Operation cost of ESF ($) | 4758.50 | 4017.93 |
| Total profit of ESF ($) | 8503.95 | 6598.34 |

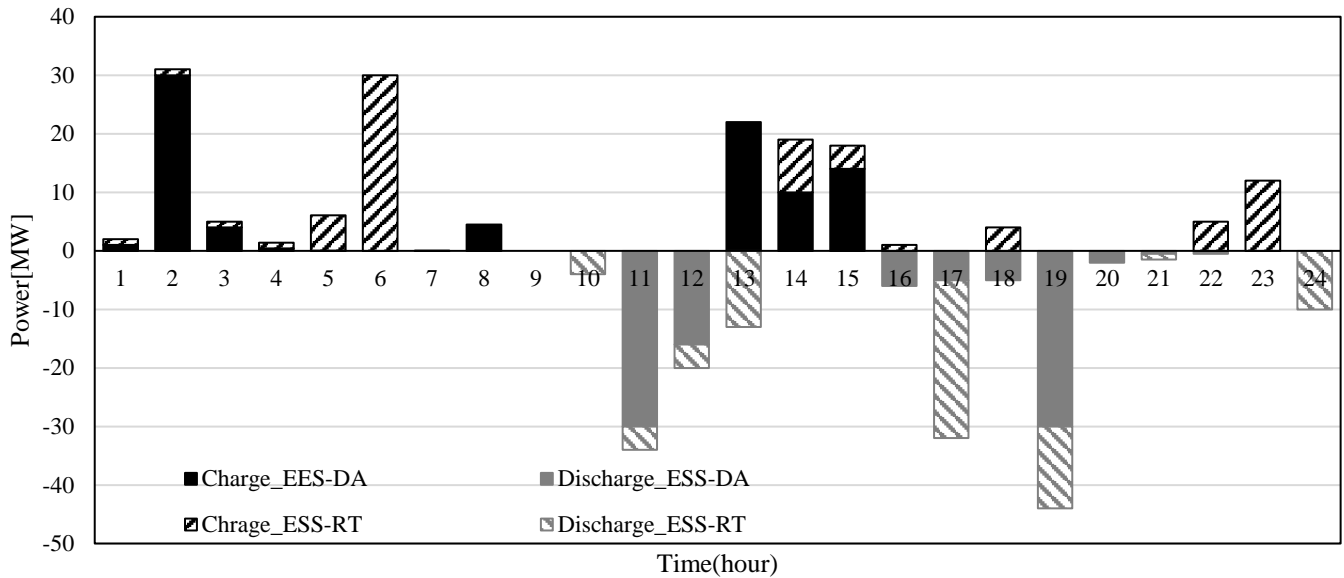

Fig. 5. Charging/discharging scheduling of ESF in power market.

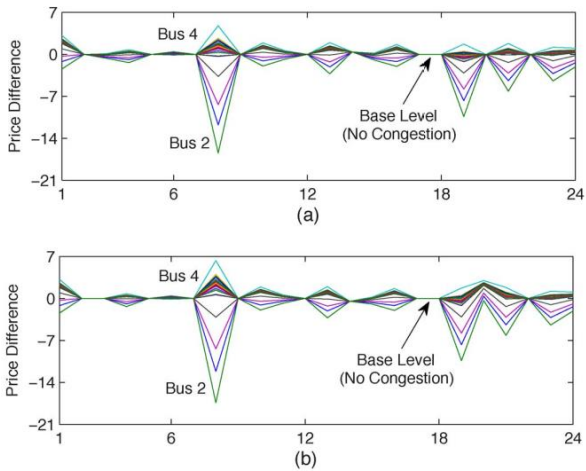

Fig. 6. LMP at 30 bus system: (a) without ESF and (b) with ESF.

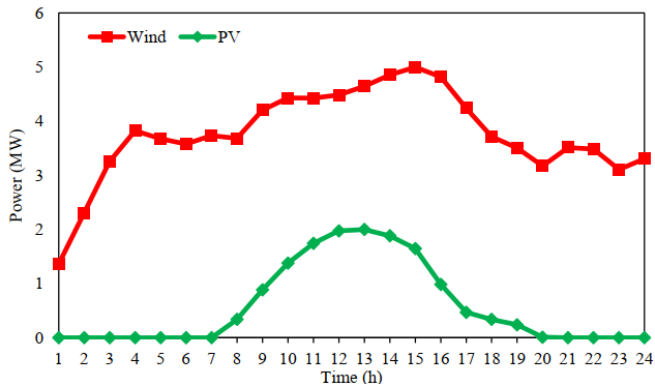

Fig. 7. Hourly predicted renewable power.

Further, to handle the uncertain parameters associated with renewable generations and the power market, a robust IGDT model has been applied to the optimization model. Hourly predicted wind and solar generations are shown in Fig. 7.

## VI. CONCLUSION

This paper addressed the strategic arbitrage of price-maker ESF facility participating in the oligopoly markets considering uncertainty. Keeping this in mind, an effective IGDT approach was employed to harness the ESF versus the uncertainty of prices. The operation decisions of the storage facility, i.e., its bids/offers in terms of both price and quantity, are strategically made in the upper-level problem. The market-clearing process under different operating conditions is modeled in the lower-level problems, aiming at maximizing the social welfare.

The proposed model was applied to obtain the strategic bidding of a merchant ESF facility based on real data from California electricity market. The actual supply and demand curves considering seven generators for each hour. The results showed that the bidding values are highly dependent upon the assumptions and the employed methods and one needs to approach the participation problem with care.

The future work of the authors includes investigating the impacts of adding other sources of uncertainty, including other sources of revenue for the facility (e.g., ancillary services or real-time markets) and considering the impact of the joint operation of storage facilities with wind/solar farms on optimal storage bidding problem.